# The Electric Location Routing Problem under Ambient Temperature

**Amin Aghalari**[1], **Darweesh Salamah**[1], **Carlos Marino**[2,3], **Mohammad Marufuzzaman**[1]
[1]Department of Industrial and Systems Engineering, Mississippi State University, MS, 39762
2CENTRUM Catolica Graduate Business School, Lima, Peru
3Pontifical Catholic University of Peru, Lima, Peru

**Abstract**

This study investigates how the location-routing decisions of the electric vehicle (EV) DCFC charging stations are impacted by the ambient temperature. We formulated this problem as a mixed-integer linear programming model that captures the realistic charging behavior of the DCFC's in association with the ambient temperature and their subsequent impact on the EV charging station location and routing decisions. Two innovative heuristics are proposed to solve this challenging model in a realistic test setting, namely, the two-phase Tabu Search-modified Clarke and Wright algorithm and the Sweep-based Iterative Greedy Adaptive Large Neighborhood algorithm. We use Fargo city in North Dakota as a testbed to visualize and validate the algorithm performances. The results clearly indicate that the EV DCFC charging station location decisions are highly sensitive to the ambient temperature, the charging time, and the initial state of charge.
*Keywords*: Electric vehicles; ambient temperature; location-routing; heuristics.

## 1 Introduction

In recent years, *electrical vehicles* (EV) have become an essential part of the manufacturing sector as the global day-by-day forced to future less dependent on nonrenewable fuel sources [1]. Sustainable transportation requires multiple efforts from different stakeholders (e.g., governments, car manufacturers, environmental advocates, and customers) to reduce the consumption of nonrenewable resources (e.g., oil, coal, and gas). EV owners will increase to around 126 million in 2030 globally and 18.7 million in the USA [1]. As more EVs take to the road, the charging station system needs to be expanded accordingly. The large-scale adoption of EVs cannot be fully realized without the adequate deployment of publicly accessible charging stations. This study is an extension of the location-routing problem to the EV area, which accounts for realistic features such as the impact of ambient temperature in the EVs' recharging process. Yang and Sun [2] proposed a mathematical model for EV battery-swap stations-based LRP such that the EVs could revisit the same swapping stations multiple times. Schiffer and Walther [3] proposed an EV LRP formulation with time windows capable of addressing a whole range of recharging options such as charging at customer sites and unique vertices, partial and full recharging. Most recently, Zhang et al. [4] proposed a hybrid heuristic algorithm that combines the binary particle swarm optimization with the variable neighborhood search to solve an EV LRP under stochastic customer demand. A comprehensive review of the LRP can be found in [5]. Despite these notable developments, past studies (e.g., [6]), especially the studies that modeled the DCFC LRPs, ignore climate variability on location-routing decisions. A recent study by Motoaki et al. [7] showed that the ambient temperature could heavily impact the DCFC charging rate. The authors stressed that considering the ambient temperature in designing the EV DCFC infrastructure in large countries like the US, where the regional climate varies significantly, could not be neglected. Unfortunately, most of the past studies (e.g., [6]) assume that the charging process of the EVs, i.e., charging rate, as a constant factor in their formulation; hence, the obtained results for the EV logistics might be altered. To fill this gap in the literature, this study extends the traditional LRPs to account for the impact of ambient temperature on the DCFC infrastructure deployment and the associated EV routing decisions.

## 2 Mathematical Model Formulation

In this section, we first proceed with the basic mathematical model formulation, referred to as **[EV]**, and then proceed to describe the model extension, referred to as **[EV-L]**. The EV basic mathematical formulation (**[EV]**) is introduced as a mixed-integer linear programming (MILP) model. We assume that there exists a linear relationship between the (*i*) travel distance and the energy consumption and (*ii*) recharging time with the amount of energy recharged. Below is a summary of the sets, parameters, and decision variables of the optimization model.

**Sets:**
- $I$: set of customers, indexed by $i \in I$
- $J$: set of potential charging station locations, indexed by $j \in J$
- $E$: set of electrical vehicles, indexed by $e \in E$
- $\{o, o'\}$: single depot and its copy
- $N$: set of all nodes, indexed by $n \in N$, where $N = \{I\} \cup \{J\} \cup \{o, o'\}$

**Parameters:**
- $f_j$: cost of installing a new charging station $j \in J$
- $d_{ij}$: distance between node $i \in I$ and $j \in J$
- $c_{ije}$: shipping cost per unit of distance between nodes $i \in I$ and $j \in J$ and EV $e \in E$
- $w_i$: demand weight for costumer $i \in I$
- $k_e$: weight capacity of EV $e \in E$
- $soc$: state of charge (SOC, in %) of an EV after getting charged in a DCFC station
- $soc_0$: initial SOC (%) of an EV at depot
- $\phi$: conversion rate of vehicle $e \in E$ which is utilized to convert the state of charge to the respective
- maximum driving distance
- $M$: a big number
- $d^{max}$: the upper bound of driving distance once EVs are fully charged

**Decision Variables:**
- $X_j$: 1 if a charging station is built in $j \in J$; 0 otherwise
- $Y_{ije}$: 1 if EV $e \in E$ traverses from node $i \in N$ to $j \in N$; 0 otherwise
- $R_{ije}$: remaining weight capacity of $e \in E$ when it arrives node $i \in N$ after leaving node $j \in N$
- $B^1_{ne}$: the maximum distance that the remaining battery power allows when EV $e \in E$ arrives at node $n \in N$
- $B^2_{ne}$: the maximum distance that the remaining battery power allows when EV $e \in E$ leaves node $n \in N$

The mathematical formulation is detailed as follows:

$$[EV] \; Minimize \; \sum_{j \in J} f_j X_j + \sum_{i \in N} \sum_{j \in N} \sum_{e \in E} c_{ije} d_{ij} Y_{ije} \tag{1}$$

Subject to,

$$\sum_{n \in N \neq \{o'\}} \sum_{e \in E} Y_{nje} = 1 \qquad \forall i \in I \tag{2}$$

$$\sum_{n \in N \setminus \{o'\}, n \neq j} \sum_{e \in E} Y_{nje} \leq M X_j \qquad \forall j \in J \tag{3}$$

$$\sum_{n' \in N \setminus \{o\}, n' \neq n} Y_{nn'e} - \sum_{n' \in N \setminus \{o'\}, n' \neq n} Y_{n'ne} = 0 \qquad \forall n \in N \setminus \{o, o'\}, e \in E \tag{4}$$

$$\sum_{n \in N \setminus \{o\}} Y_{one} - \sum_{n \in N \setminus \{o'\}} Y_{noe} = 0 \qquad \forall e \in E \tag{5}$$

$$\sum_{n \in N \setminus \{o\}} Y_{one} \leq 1 \qquad \forall e \in E \qquad (6)$$

$$\sum_{n \in N \setminus \{o'\}, n \neq j} R_{nje} = \sum_{n \in N \setminus \{o\}, n \neq j} R_{jne} \qquad \forall j \in J, e \in E \qquad (7)$$

$$\sum_{n \in N \setminus \{o \cup i\}} R_{ine} \leq \sum_{n \in N \setminus \{o' \cup i\}} R_{ine} - w_i \sum_{n \in N \setminus \{o' \cup i\}} Y_{nie} + k_e (1 - \sum_{n \in N \setminus \{o' \cup i\}} Y_{nie}) \qquad \forall i \in I, e \in E \qquad (8)$$

$$0 \leq R_{nn'e} \leq k_e Y_{nn'e} \qquad \forall n \in N \setminus \{o'\}, n' \in N \setminus \{o\}, n \neq n', e \in E \qquad (9)$$

$$B^1_{n'e} + d_{nn'} Y_{nn'e} \leq B^2_{ne} + d^{max}(1 - Y_{nn'e}) \qquad \forall n \in N \setminus \{o'\}, n' \in N \setminus \{o\}, n \neq n', e \in E \qquad (10)$$

$$B^1_{n'e} + d_{nn'} Y_{nn'e} \geq B^2_{ne} + d^{max}(1 - Y_{nn'e}) \qquad \forall n \in N \setminus \{o'\}, n' \in N \setminus \{o\}, n \neq n', e \in E \qquad (11)$$

$$B^2_{oe} = \phi soc_0 \qquad \forall e \in E \qquad (12)$$

$$B^2_{oe} = \phi soc X_j \qquad \forall j \in J, e \in E \qquad (13)$$

$$B^1_{ie} = B^2_{ie} \qquad \forall i \in I, e \in E \qquad (14)$$

$$R_{ije}, B^1_{ne}, B^2_{ne} \geq 0 \qquad \forall (i,j) \in N, n \in N, e \in E \qquad (15)$$

$$X_j, Y_{ije} \in \{0,1\} \qquad \forall i \in N, j \in J\ e \in E \qquad (16)$$

The objective function of **[EV]** minimizes the total cost associated with opening EV charging stations and the driving distance costs within a planning horizon. Constraints (2) ensure that each customer site $i \in I$ is visited by exactly one EV. Constraints (3) ensure that EVs could get recharged at a specific charging station $j \in J$ if only it is located. Constraints (4) enforce the flow balance for each EV's e\inE in the customer sites and the charging stations. Constraints (5) guarantee that a utilized EV $e \in E$ should return to the depot at the end of the respective trip. Constraints (6) limit the number of trips that an EV $e \in E$ can start from the depot. Constraints (7) ensure that at any charging station $j \in J$, the remaining weight capacity of the EVs does not change ($w_j = 0; \forall j \in J$). Constraints (8) update the remaining weight capacity of the EVs based on the nodes visited. Constraints (9) enforce that the remaining weight capacity of the EVs is less than the EV maximum weight capacity and also be greater than zero in all the visited nodes by the EVs. Constraints (10) and (11) update the battery power level of the EVs based on the nodes visited. Constraints (12) and (13) detail the SOC when the EV $e \in E$ starts its trip from the depot and when it visits a charging station. Constraints (14) ensure that the battery level of the EVs $e \in E$ remains unchanged when they visit a customer node $i \in I$ in the network. Constraints (15) and (16) enforce nonnegativity and binary restrictions for the decision variables.

Model **[EV]** assumes that the SOC of a fast charger drops linearly. However, the actual fast charging process is non-linear and is a function of initial SOC and ambient temperature [7]. The simplified linearized SOC assumption may provide an overestimated duration for the DCFC's. This sub-section introduces model **[EV-L]** by alleviating this drawback from model [EV]. Let us define $s\hat{o}c(c,t)$ to predict the SOC of an EV $e \in E$, which is a function of charging time $t$ and ambient temperature $c$ (in Celsius). We further define $\lambda_0, \lambda_1$, and $\lambda_2$ to be the coefficient estimates, and $soc_i$ the initial value of the SOC of an EV. Inspired from the study of [7], the following SOC estimation is provided.

$$s\hat{o}c(c,t) = e^{\lambda_2 t} soc_i + \left(\frac{\lambda_0 + \lambda_1}{\lambda_2}\right) = \mu_1 soc_i + \mu_2$$

Using the conversion rate ($\phi$), we know that $soc_i = \frac{B^1_{je}}{\phi}$; hence, we obtain the following:

$$B^2_{je} \leq soc \phi X_j \rightarrow B^2_{je} \leq \left(\mu_1 \frac{B^1_{je}}{\phi} + \mu_2\right) \phi X_j \rightarrow B^2_{je} \leq \mu_1 B^1_{je} X_j + \mu_2 \phi X_j \qquad \forall j \in J, e \in E \qquad (17)$$

Constraints (17) are nonlinear. We introduce a new variable $\{Z_{je} | \forall j \in J, e \in E\}$ to replace the $B^1_{je} X_j$ term.

$$Z_{je} \leq d^{max} X\_j \qquad \forall j \in J, e \in E \qquad (18)$$
$$Z_{je} \leq B_{je}^1 \qquad \forall j \in J, e \in E \qquad (19)$$
$$Z_{je} \geq B_{je}^1 - d^{max}(1 - X_j) \qquad \forall j \in J, e \in E \qquad (20)$$
$$Z_{je} \geq 0 \qquad \forall j \in J, e \in E \qquad (21)$$

With this, model **[EV]** can be extended as follows, referred to as **[EV-L]**:

$$[EV - l] \text{ Minimize} \sum_{j \in J} f_j X_j + \sum_{i \in N} \sum_{j \in N} \sum_{e \in E} c_{ije} d_{ij} Y_{ije} \qquad (22)$$

subject to: (2)-(12), (14)-(16), and (17)-(21).

## 3 Solution Methodology

Both basic (**[EV]**) and extended (**[EV-L]**) formulations developed in this study are indeed variants of the classical location-routing problems; As such, they can be considered as an *Np*-hard problem. Our initial experimentation with the GUROBI solver exposes its inability to solve the largest instances of problem **[EV-L]** in a reasonable timeframe. This study proposes two heuristic techniques, namely, the two-phase Tabu Search-modified Clarke and Wright Savings heuristic (TS-MCWS) and the Sweep-based Iterated Greedy Adaptive Large Neighborhood algorithm (SIGALNS), to solve model **[EV-L]** efficiently.

### 3.1 Framework of TS-MCWS and SIGALNS Algorithm

The TS-MCWS heuristic, which combines the Tabu Search (TS) algorithm with a modified version of the Clarke and Wright Savings method. Within this two-phase algorithm, the TS algorithm is used to determine the location of the charging stations and then given the selected charging stations, the modified Clarke and Wright Savings (MCWS) method is used to finding the routing decisions. Two algorithms collaborate iteratively to provide an efficient solution for the model **[EV-L]** as follows. Note that, $S_0$, $S$, $S^*$ stand for an initial solution, the current solution, and the best-known solution, respectively.

- **Step 1**: The initial number of the charging stations, $N_j$, is set to one.
- **Step 2**: Using the radius covering algorithm, $N_j$ charging stations are selected. Then, using the MCWS procedure and considering the selected charging stations, the routing plan, $S_0$, to satisfy the customer demands, are obtained. By doing so, the current solution $S$ and the best-known solution $S^*$ are set to the initial solution $S_0$.
- **Step 3**: By applying the TS procedure on $S$, the current solution is updated. If $Z(S) < Z(S^*)$, the best-known solution is updated, $S^* \leftarrow S$.
- **Step 4**: If a pre-specified number of iterations without improvement in the objective function value of the best-known solution $N_{itr}$ has reached or all the charging stations have located, the TS-MCWS heuristic is terminated. Otherwise, set $N_j \leftarrow N_j + 1$ and proceed to **Step 2**. In our experiments, if the
number of costumers $|I| \leq 75$, we set $N_{itr} = 5$; otherwise, we set $N_{itr} = 10$.

The SIGALNS heuristic composes of three components, namely, the modified Sweep heuristic, the Iterated Greedy (IG), and the Adaptive Large Neighborhood Search (ALNS) algorithm, to solve model **[EV-L]**. Within this algorithm, the IG algorithm is used to determine the location of the charging stations and then given the selected charging stations, the ALNS algorithm is used to finding the routing decisions (see **Algorithm 1**).

## 4 Computational performance and Case Study

The efficiency of the proposed algorithms in solving model **[EV-L]** are evaluated on different test instances. To do so, a new set of problem instances with various sizes, in terms of the number of the customers and the available EVs, are generated. In total 150 locations with a considerable population around center of the test region (Fargo) are considered as the potential locations of the customers. Next, varying the size of the customers, 10 different test instances are generated. The computational performance of the proposed algorithms under these generated test instances are discussed below.

**Algorithm 1: The framework of SIGALNS algorithm**

**Input:** The distance between all the nodes $d_{n,n'}$, $\forall (n, n') \in N$, customers' demand $w_i$, the wight capacity of EVs $k_e$, $\forall e \in E$, battery driving range of EVs $s\hat{o}c(c, t)$, where $c$ is charging time and $t$ is the ambient temperature

Implement modified sweep algorithm to obtain a initial solution $S_0$
$S \leftarrow S_0, S^* \leftarrow S_0$
Starting the initial value of the removal and insertion operators for the ALNS algorithm
$iter \leftarrow 1$
**while** $iter < ITR^{SIGALNS}$ **do**
  $S' \leftarrow S$
  Remove all the located charging stations from $S'$
  Apply the iterated greedy algorithm to $S'$ to find out the updated located charging stations
  Apply the ALNS algorithm to $S'$ to update the routing plans
  **if** the acceptance criterion is satisfied **then**
   | $S \leftarrow S'$
  **end**
  **if** $Z(S) < Z(S^*)$ **then**
   | $S^* \leftarrow S$
  **end**
  $iter \leftarrow iter + 1$
**end**
**Output:**
  The best-known solution: $S^*$

The first two columns in Tables 1-3 represent the problem instances and the number of respective customers. Next, $|J_l|$, $Best$, $Average$, and $T$ to represent the number of charging stations opened, the best feasible solution, the average feasible solution, and the average running time (in seconds) of the investigated algorithms, respectively. The last column, $gap(\%)$, in Tables 1-3 represents the difference between the best feasible solution found by the two heuristics and is computed as follows: ($Best_2$ $Best_1$)/$Best_1$. It is worth mentioning that we run each test instance *five* times to obtain the average solution and running time reported in Tables 1-3. Results indicate that SIGALNS algorithm provides a high-quality feasible solution over the TS-MCWS algorithm. Apart from the improvement in solution quality in the SIGALNS algorithm, the best-known solution is achieved faster than the TS-MCWS algorithm.

Table 1: Performance of SIGALNS and TS-MCWS when temperature is $-10\,°C$

|  |  | SIGALNS | | | | | TS-MCWS | | | |
| --- | --- | --- | --- | --- | --- | --- | --- | --- | --- | --- |
|  | $|I|$ | $|J_l|$ | $Best_1$ | $Average_1$ | $T1(s)$ | $|J_l|$ | $Best_2$ | $Average_2$ | $T2(s)$ | $gap(\%)$ |
| Average over 10 instances | 57.5 | 11 | 32,773.6 | 32,860.8 | 49.9 | 11.4 | 33,907.9 | 34,035.3 | 132.7 | 3.1 |

Table 2: Performance of SIGALNS and TS-MCWS when temperature is $10\,°C$

|  |  | SIGALNS | | | | | TS-MCWS | | | |
| --- | --- | --- | --- | --- | --- | --- | --- | --- | --- | --- |
|  | $|I|$ | $|J_l|$ | $Best_1$ | $Average_1$ | $T1(s)$ | $|J_l|$ | $Best_2$ | $Average_2$ | $T2(s)$ | $gap(\%)$ |
| Average over 10 instances | 57.5 | 5 | 16,866.6 | 16,925.5 | 43.9 | 5.5 | 17,751.2 | 17,850.7 | 116.9 | 3.8 |

Table 3: Performance of SIGALNS and TS-MCWS when temperature is $30\,°C$

|  |  | SIGALNS | | | | | TS-MCWS | | | |
| --- | --- | --- | --- | --- | --- | --- | --- | --- | --- | --- |
|  | $|I|$ | $|J_l|$ | $Best_1$ | $Average_1$ | $T1(s)$ | $|J_l|$ | $Best_2$ | $Average_2$ | $T2(s)$ | $gap(\%)$ |
| Average over 10 instances | 57.5 | 1.7 | 8,055.8 | 8,108.6 | 28 | 1.8 | 8,460.6 | 8,636.7 | 67 | 3.5 |

Finally, we perform a set of sensitivity analysis to evaluate the performance of the proposed model. From Figure 1(a), it can be observed that the EV DCFC charging station location decisions are sensitive to the ambient temperature and charging time. For instance, when the charging time of the EVs is set to its base value (80 minutes) and the ambient temperature decreases from -10 °C to 10 °C, the selection of the charging stations increases by approximately 141%. Figure 1(b) indicates the impact of variation in the ambient temperature on the overall system cost. As can be seen from the figure, increasing the ambient temperature, the median value of the overall system cost decreases. For instance, when the ambient temperature increases from $-10$ °C to 10 °C, the overall system cost drops by approximately 23.1%. The results clearly indicate that the EV DCFC charging station location decisions are highly sensitive to the ambient temperature and charging time.

# 5 Conclusion and Future Research Directions

This study proposes a mixed-integer linear programming model to minimize the EV DCFC infrastructure and the associated routing decisions under fluctuating ambient temperature. Two highly customize heuristic approaches are proposed to efficiently solve the optimization model. This research can be extended in several directions. First, it would be interesting to see how the stochasticity associated with different input parameters (e.g., charging rate, customer demand) impact the EV DCFC LRP. Next, efforts will continue to develop rigorous techniques such as decomposition methods [8–10] to improve the quality of the solutions.

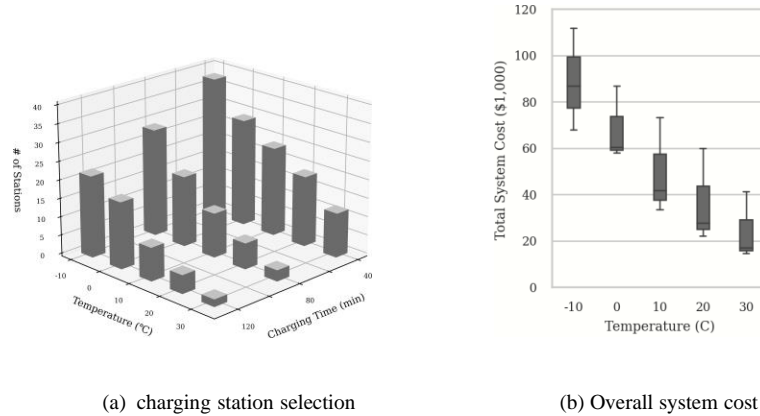

(a) charging station selection  (b) Overall system cost

Figure 1: Impact of temperature and charging time on EV DCFC for largest instance with 150 customers.

## References

[1] Matulka, R., 2014. Timeline: history of the electric car. *Retrieved from United States Department of Energy Website: http://energy. gov/articles/historyelectric-car*.

[2] Yang, J. and Sun, H., 2015. Battery swap station location-routing problem with capacitated electric vehicles. Computers & Operations Research, 55, pp.217-232.

[3] Schiffer, M. and Walther, G., 2017. The electric location routing problem with time windows and partial recharging. European Journal of Operational Research, 260(3), pp.995-1013.

[4] Zhang, S., Chen, M. and Zhang, W., 2019. A novel location-routing problem in electric vehicle transportation with stochastic demands. Journal of Cleaner Production, 221, pp.567-581.

[5] Prodhon, C. and Prins, C., 2014. A survey of recent research on location-routing problems. European Journal of Operational Research, 238(1), pp.1-17.

[6] Kabli, M., Quddus, M.A., Nurre, S.G., Marufuzzaman, M. and Usher, J.M., 2020. A stochastic programming approach for electric vehicle charging station expansion plans. International Journal of Production Economics, 220, p.107461.

[7] Motoaki, Y., Yi, W. and Salisbury, S., 2018. Empirical analysis of electric vehicle fast charging under cold temperatures. Energy Policy, 122, pp.162-168.

[8] Aghalari, A., Nur, F. and Marufuzzaman, M., 2020. A Bender's based nested decomposition algorithm to solve a stochastic inland waterway port management problem considering perishable product. International Journal of Production Economics, 229, p.107863.

[9] Aghalari, A., Nur, F. and Marufuzzaman, M., 2020. Solving a stochastic inland waterway port management problem using a parallelized hybrid decomposition algorithm. Omega, p.102316.

[10] Marufuzamman, M., Aghalari, A., Buchanan, R.K., Rinaudo, C.H., Houte, K.M. and Ranta, J.H., 2020. Optimal Placement of Detectors to Minimize Casualties in an Intentional Attack. IEEE Transactions on Engineering Management.